\documentclass[12pt]{amsart}

\usepackage{amssymb}

\usepackage{enumitem}

\usepackage{graphicx}

\usepackage{mathrsfs}
\usepackage{mathtools}
\usepackage{graphicx,scalerel}

\makeatletter
\@namedef{subjclassname@2020}{%
  \textup{2020} Mathematics Subject Classification}
\makeatother

\newtheorem{theorem}{Theorem}[section]

\newtheorem{proposition}[theorem]{Proposition}

\theoremstyle{definition}
\newtheorem{definition}[theorem]{Definition}
\newtheorem{remark}[theorem]{Remark}
\newtheorem{example}[theorem]{Example}

\numberwithin{equation}{section}

\newcommand\wh[1]{\hstretch{2}{\hat{\hstretch{.5}{#1}}}}

\begin{document}

%%%%% To ease editing, for IMPAN journals add:

\baselineskip=17pt

%%%%%%%%%%%%%%%%

\title[Real Structures]{Real Structures. An introduction to a general approach}

\author[T. Valent]{Tullio Valent}
\address{Department of Mathematics ``Tullio Levi Civita''\\ University of Padua\\
Via Trieste 63\\
35121 Padua, Italy}
\email{tullio.valent@unipd.it}

%\author[J. K. Nowak]{Jan Krzysztof Nowak}
%\address{Institute of Mathematics\\ Jagiellonian University\\
%{\L}ojasiewicza 6\\
%30-348 Krak\'ow, Poland}
%\email{jk.nowak@im.uj.edu.pl}

\date{}

\begin{abstract}
In this paper we attempt to present a very general approach to the study of structures (somehow) defined on a set $X$ by a family of maps $d: X \times X \mapsto \mathbb{R}^+$.

It will be shown how the assignment of a preorder $\prec_{\Pi}$ on a set $\Pi$ of families of maps from $X \times X$ into $\mathbb{R}^+$ defines a structure on $X$.
The structures obtained in this way will be called \emph{real structures}. For real structures on two different sets we study when they are \emph{of the same type}. An answer to this question will allow to introduce the notion of \emph{morphism}, and then to give the definition of the initial real structure with respect to a family of maps, and so also the definition of product real structure.

Various examples of preorders, and hence of real structures, will be exhibited and discussed.

A few examples of morphisms will be proposed.
\end{abstract}



\maketitle

\section{Introduction}\label{intro}
We start by considering a set $\Pi$ of families of maps $d : X \times X \mapsto \mathbb{R}^+$, where $X$ is any set. $\Pi$ may be the set of all families of maps $d : X \times X \mapsto \mathbb{R}^+$. An important choice of $\Pi$ is the family of all pseudo-metrics on $X$. It is indeed well known that all uniformities for $X$ can be defined by a family of pseudo-metrics on $X$: that is, the family of all uniformly continuous pseudo-metrics.  

Another interesting choice of $\Pi$ seems to us the family of all \emph{weak pseudo-metrics} on $X$. By  a weak pseudo-metrics on $X$ we mean a symmetric map $d : X \times X \mapsto \mathbb{R}^+$ that satisfy the (triangle) inequality 
$$ d(x_1,x_2) \leq d(x_1,x) +d(x,x_2) \,\, \text{for all} \, \, x_1,x_2,x \in X$$
and vanishes in at least a point of the diagonal $\Delta$ of $X \times X$.

We emphasize the fact that the weak pseudo-metrics are not requested to vanish on the whole of $\Delta$ (unlike in the usual meaning of pseudo-metric). It is worth recalling (cfr. \cite{val23}) that all topologies on a set $X$, all pseudo-uniformities for $X$, and many other structures can be defined by a suitable family of weak pseudo-metrics (see Sect.~\ref{ex2:preord}).  

In Sect. \ref{realstruct} we will expose the preliminaries to our definition of real structure, founded on an assigned preorder $\prec_{\Pi}$ on  $\Pi$.

In Sect. \ref{realspaces} we will deal with a definition of \emph{``type" of a real structure}: it will be necessary to decide when real structures defined on different sets are of the same type. Then we will be leaded to the notion of \emph{morphism}, which is important because it allows to define the real structure induced in a subset, the \emph{product real structure}, and, more in general, to define the \emph{initial real structure} with respect to a family of maps.

Subsequently, a few examples of preorders $\prec_{\Pi}$ and explicit examples of real structures will be presented. Moreover, we will linger over various examples of morphisms between real spaces.

\section{Real structures on a set $X$. Preliminaries.}\label{realstruct}

Let $\Pi$ be a set of families $\mathscr{P}$ of maps $d : X \times X \mapsto \mathbb{R}^+$.
On $\Pi$ let be assigned a preorder $\prec_{\Pi}$ which has the property 
\begin{equation}\label{preorder}
	\mathscr{P}' \prec_{\Pi} \mathscr{P} \Leftrightarrow \{d'\} \prec_{\Pi}  \mathscr{P} \quad \forall d' \in \mathscr{P}'.
\end{equation}
In the case  $\Pi$ is the set of all families of maps $d : X \times X \mapsto \mathbb{R}^+$ the preorder could be denoted, more simply, with the symbol
$\prec$.
As, obviously, since $\prec_{\Pi}$ is a preorder, $\mathscr{P} \prec_{\Pi} \mathscr{P}$, we  have $\{d\} \prec_{\Pi} \mathscr{P}$ for all $d \in \mathscr{P}$. Note that from \eqref{preorder} it follows that
$$ \mathscr{P}' \subseteq \mathscr{P} \Rightarrow \mathscr{P}' \prec_{\Pi} \mathscr{P}.$$
If $\mathscr{P}' \prec_{\Pi} \mathscr{P}$ we will say that $\mathscr{P}$ \emph{absorbs}  $\mathscr{P}'$. The preorder $\prec_{\Pi}$ induces on $\Pi$ an equivalence relation: $\mathscr{P}$ is equivalent to  $\mathscr{P}'$ when $\mathscr{P} \prec_{\Pi} \mathscr{P}'$
and $\mathscr{P}' \prec_{\Pi} \mathscr{P}$. For any element $\mathscr{P}$ of $\Pi$ we will set
$${\scriptstyle\prec_{\Pi}}(\mathscr{P})\coloneqq \{d : X \times X \mapsto \mathbb{R}^+ : \{d\}\prec_{\Pi} \mathscr{P} \}.$$
Obviously $\mathscr{P} \subseteq {\scriptstyle\prec_{\Pi}}(\mathscr{P})$. Note that ${\scriptstyle\prec_{\Pi}}(\mathscr{P})$ \emph{is the greatest sub-family of $\Pi$ which is equivalent to $\mathscr{P}$}; therefore it can be identified with the equivalence class containing $\mathscr{P}$. It will be called the \emph{real structure generated by $\mathscr{P}$} (through the preorder$\prec_{\Pi}$).

\begin{theorem}
	A sub-family $\mathscr{R}$ of $\Pi$ is a real structure generated (through a preorder $\prec_{\Pi}$) by some family $\mathscr{P}$ of $\Pi$ if and only if
	$$ (R) \qquad   {\scriptstyle\prec_{\Pi}}(\mathscr{R}) \subseteq \mathscr{R} \qquad [\text{which is equivalent to }  {\scriptstyle\prec_{\Pi}}(\mathscr{R})=\mathscr{R} ],$$
	namely $\{d\} \prec_{\Pi}\mathscr{R} \Rightarrow d \in \mathscr{R}$.
\end{theorem}
\begin{proof}
	Observe that ${\scriptstyle\prec_{\Pi}}(\mathscr{P}) \subseteq {\scriptstyle\prec_{\Pi}}({\scriptstyle\prec_{\Pi}}(\mathscr{P}))$ and that $${\scriptstyle\prec_{\Pi}}({\scriptstyle\prec_{\Pi}}(\mathscr{P})) \prec_{\Pi} {\scriptstyle\prec_{\Pi}}(\mathscr{P})\prec_{\Pi} \mathscr{P},$$ which implies ${\scriptstyle\prec_{\Pi}}({\scriptstyle\prec_{\Pi}}(\mathscr{P})) \subseteq {\scriptstyle\prec_{\Pi}}(\mathscr{P})$.
	
	Then if $\mathscr{R} = {\scriptstyle\prec_{\Pi}}(\mathscr{P})$ we have ${\scriptstyle\prec_{\Pi}}(\mathscr{R}) \subseteq \mathscr{R}$.
	
	Conversely, if ${\scriptstyle\prec_{\Pi}}(\mathscr{R}) \subseteq \mathscr{R}$ then $\mathscr{R}$ is generated by itself.    
\end{proof}
The identification of the equivalence class containing $\mathscr{R}$ with ${\scriptstyle\prec_{\Pi}}(\mathscr{R})$ leads us to the following

\begin{definition}
	A \emph{real structure on $X$ of type $\prec_{\Pi}$} is a sub-family $\mathscr{R}$ of $\Pi$ having the property $(R)$. The pair $(X,\mathscr{R})$ will be called a \emph{real space of type $\prec_{\Pi}$}. 
\end{definition}

\section{Real spaces of the same type. Morphisms.}\label{realspaces}
\begin{definition}
	Two real spaces $(X,\mathscr{R}_X)$ and $(Y,\mathscr{R}_Y)$ are \emph{of the same type} if there are two maps $f: X \mapsto Y$ and $g: Y \mapsto X$ such that 
	$\mathscr{R}_Y \circ \overset{(2)}{f}  \subseteq \mathscr{R}_X$ and $\mathscr{R}_X \circ \overset{(2)}{g}  \subseteq \mathscr{R}_Y$, where
	$$\mathscr{R}_Y \circ \overset{(2)}{f} \coloneqq \{d_Y \circ \overset{(2)}{f} : d_Y \in \mathscr{R}_Y \}, \quad 
	\mathscr{R}_X \circ \overset{(2)}{g} \coloneqq \{d_X \circ \overset{(2)}{g} : d_X \in \mathscr{R}_X \}.$$
\end{definition}

A map $f: X \mapsto Y$ such that 	$\mathscr{R}_Y \circ \overset{(2)}{f}  \subseteq \mathscr{R}_X$ will be called a \emph{morphism of $X$ into $Y$} and a map 
$g: Y \mapsto X$ such that $\mathscr{R}_X \circ \overset{(2)}{g}\subseteq \mathscr{R}_Y$ will be called a \emph{morphism of $Y$ into $X$} (with respect to the real structures $\mathscr{R}_X$ and $\mathscr{R}_Y$).

\begin{remark}
	If $\mathscr{P}_X$ generates $\mathscr{R}_X$ and $\mathscr{P}_Y$ generates $\mathscr{R}_Y$ then a map $f : X \mapsto Y$ is a morphism if and only if
	\begin{equation}\label{morphgen}
		d \circ \overset{(2)}{f} \prec_{\Pi_X} \mathscr{P}_X \qquad \forall \, \{d\} \prec_{\Pi_Y} \mathscr{P}_Y \qquad (\text{i.e} \,\, \mathscr{P}_Y \circ \overset{(2)}{f}\prec_{\Pi_X} \mathscr{P}_X),
	\end{equation} 
	where $\prec_{\Pi_X}$ and $\prec_{\Pi_Y}$ are the preorders on $X$ and $Y$.
\end{remark}
\begin{proof}
	Since $\mathscr{R}_Y = \{d: Y \times Y \mapsto \mathbb{R}^+ : \{d\} \prec_{\Pi_Y} \mathscr{P}_Y\}$ we have $\mathscr{R}_Y \circ \overset{(2)}{f} = 
	\{d \circ  \overset{(2)}{f} :\{d\} \prec_{\Pi_Y} \mathscr{P}_Y  \}$; thus \eqref{morphgen} implies $\mathscr{R}_Y \circ \overset{(2)}{f} \prec_{\Pi_X} \mathscr{P}_X$ and so  $\mathscr{R}_Y \circ \overset{(2)}{f} \subseteq \mathscr{R}_X$. Conversely, $\mathscr{R}_Y \circ \overset{(2)}{f} \subseteq \mathscr{R}_X$ implies  $\mathscr{P}_Y \circ \overset{(2)}{f} \subseteq \mathscr{P}_X$, and this implies   $\mathscr{P}_Y \circ \overset{(2)}{f} \prec_{\Pi_X} \mathscr{P}_X$.
\end{proof}
\begin{definition}
	Let $X$ be a set, $(Y,\mathscr{R}_Y)$  a real space and $f:X \mapsto Y$ a map. The real structure on $X$, of the same type of $\mathscr{R}_Y$, generated by the family $\mathscr{R}_Y \circ \overset{(2)}{f}$ is (clearly) the smallest one for which $f$ is a morphism. We will say that it is the \emph{initial real structure on $X$ with respect to $f$}. More in general, if $X$ is a set, $(Y_i,\mathscr{R}_{Y_i})_{i \in I}$ a family of real spaces of the same type, and, for each $i \in I$, $f_i : X \mapsto Y_i$ a map, the \emph{initial real structure on $X$ with respect to the family $(f_i)_{i \in I}$} [i.e. the smallest real structure on X, of the same type of the real structures $\mathscr{R}_{Y_i}$, for which each $f_i$ is a morphism] is generated by the family
	$(\mathscr{R}_{Y_i} \circ \overset{(2)}{f_i})_{i \in I}$.
	
	Thus, we can define the notions of \emph{real structure induced on a subset} (as the initial real structure with respect to the inclusion map) and the \emph{product real structure} (as the initial real structure with respect to the projections).
\end{definition}

\section{$\Delta$-local filters on $X \times X$, and topologies induced on $X$.}\label{deltalocal}

\begin{definition}
	Let $X$ be a set and, for each $x \in X$, $\mathscr{F}_x$ be a filter on $X \times X$ whose elements contain the point $(x,x)$. The union
	$$\mathscr{F} \coloneqq \bigcup_{x \in X} \mathscr{F}_x$$
	will be called a \emph{$\Delta$-local filter on $X\times X$}.
\end{definition}

Evidently we have 
$$ \mathscr{F}_x =\{U \in \mathscr{F} : (x,x) \in U\}.$$

$\mathscr{F}$ is an union of filters, but (in general) is not a filter. Observe that $\mathscr{F}$ is a filter on $X \times X$ when every element of $\mathscr{F}$ contains $\Delta$.

\begin{definition}
	Any family $\mathscr{P}$ of maps $d : X \times X \mapsto \mathbb{R}^+$ defines the $\Delta$-local filter $\wh{\mathscr{F}}(\mathscr{P})$, where
	$\wh{\mathscr{F}}(\mathscr{P}) \coloneqq \bigcup_{x \in X} \wh{\mathscr{F}_x}(\mathscr{P})$, with $\wh{\mathscr{F}_x}(\mathscr{P})$  the filter on $X \times X$ generated by  $\{d^{\leftarrow}\left([0,\varepsilon[\right): d \in \mathscr{P}, \varepsilon > d(x,x)\}$.  
\end{definition}  
\begin{proposition}
	Any family $\mathscr{P}$	of maps  $d : X \times X \mapsto \mathbb{R}^+$ defines a topology 	$\wh{\mathscr{\tau}}(\mathscr{P})$ in the following way. For each $x \in X$, $d \in \mathscr{P}$ and $\varepsilon > d(x,x)$ we set $U_{d,\varepsilon}(x)=\{\xi \in X: d(\xi,x) < \varepsilon, d(x,\xi)< \varepsilon)\}$.
	The family $\{U_{d,\varepsilon}(x) : d \in \mathscr{P},\varepsilon > d(x,x)\}$ is a pre-base for a filter, say $\mathscr{F}(x)$, on $X$. Then we set 
	$\wh{\mathscr{\tau}}(\mathscr{P})\coloneqq \{A \subseteq X: A \in \mathscr{F}(a) \, \forall a \in A\}$. $\wh{\mathscr{\tau}}(\mathscr{P})$ is a topology on X.  
\end{proposition}

\begin{proposition}
	Any $\Delta$-local filter $\mathscr{F}$ on  $X \times X$ induces a topology $\wh{\mathscr{\tau}}(\mathscr{F})$ on $X$.
\end{proposition}

\begin{proof}
	Consider, for every $x \in X$, the set $\mathscr{F}[x] \coloneqq \{U[x]: u \in \mathscr{F}_x\}$, where $U[x] \coloneqq \{\xi \in X : (\xi,x),(x,\xi) \in U\}$.
	Of course, $x \in U[x]$. Now set
	$$\wh{\mathscr{\tau}}(\mathscr{F}) \coloneqq \{A \subseteq X : A \in \mathscr{F}[a] \,\, \forall a \in A\}.$$
	Observe that $\wh{\mathscr{\tau}}(\mathscr{F})$ is closed under finite intersections, because if
	
	\begin{equation}\label{int:1}
		A_1 \in \mathscr{F}[a_1] \,\, \forall a_1 \in A_1, \quad \text{and} \quad A_2 \in \mathscr{F}[a_2] \,\, \forall a_2 \in A_2,
	\end{equation}
	
	then
	
	\begin{equation}\label{int:2}
		A_1 \cap A_2 \in \mathscr{F}[a] \,\, \forall a \in A_1 \cap A_2.
	\end{equation}
	
	Indeed, from \eqref{int:1} it follows that there are $U_1, U_2 \in \mathscr{F}$ such that 
	$$A_1 \in U_1[a_1] \, \, \forall a_1 \in A_1,\quad \text{and} \quad A_2 \in U_2[a_2] \, \, \forall a_2 \in A_2,$$ which implies $A_1 \cap A_2 \in U_1[a_1] \cap U_2[a_2] \, \, \forall (a_1,a_2) \in A_1 \times A_2$, and so, putting $U = U_1 \cap U_2$, it is easily seen that \eqref{int:2} is true. Of course, $U$ may be empty; in this case $A_1 \cap A_2 = \emptyset$ and hence \eqref{int:2} holds trivially.
	
	Since, obviously, $\wh{\mathscr{\tau}}(\mathscr{F})$ is closed under every union, we can conclude that $\wh{\mathscr{\tau}}(\mathscr{F})$ is a topology on $X$, which we refer to as the topology induced by $\mathscr{F}$. 
\end{proof}

\section{A few examples of preorders $\prec_{\Pi}$.}\label{ex:preord}
Let $X$ be any set, and let $\mathscr{P},\mathscr{P}'$ be two families of $\mathbb{R}^+$-maps defined in $X \times X$. As done in Sect.~\ref{deltalocal}, we denote by $\wh{\mathscr{\tau}}(\mathscr{P})$ and $\wh{\mathscr{\tau}}(\mathscr{P}')$ the topologies on $X$ defined by $\mathscr{P}$ and $\mathscr{P}'$, while 
$\wh{\mathscr{F}}(\mathscr{P})$ and $\wh{\mathscr{F}}(\mathscr{P}')$ will be the $\Delta$-local filters on $X \times X$ defined by $\mathscr{P}$ and $\mathscr{P}'$.
We recall that
$$\wh{\mathscr{F}}(\mathscr{P}) = \bigcup_{x \in X} \wh{\mathscr{F}_x}(\mathscr{P}),$$
where $\wh{\mathscr{F}_x}(\mathscr{P})$ is the filter generated by the family of sets $\{(x_1,x_2) \in X \times X : d(x_1,x_2)<\varepsilon\}$, with $d \in \mathscr{P}$ and $\varepsilon > d(x,x)$.

\begin{example}\label{ex:top}
	We say that $\mathscr{P}' \prec_{\Pi} \mathscr{P}$ \emph{in the topological way} if $\wh{\mathscr{\tau}}(\mathscr{P}') \subseteq \wh{\mathscr{\tau}}(\mathscr{P})$, and so $\mathscr{P}$ and $\mathscr{P}'$ are equivalent if they define the same topology on $X$. Thus, a topological structure $\mathscr{R}_{\tau}$, thought as a ``real" structure, is (in the case where $\Pi$ is the set of all maps $d : X \times X \mapsto \mathbb{R}^+$) the greatest family of maps $d : X \times X \mapsto \mathbb{R}^+$ defining a topology on $X$. It is easy to prove that any topological ``real" structure $\mathscr{R}_{\tau}$ has the following properties:
	\begin{equation}
		\begin{cases*}
			d_1 \leq d, d \in \mathscr{R}_{\tau} \Rightarrow d_1 \in  \mathscr{R}_{\tau}\\
			d \in  \mathscr{R}_{\tau} \Rightarrow \alpha d \in  \mathscr{R}_{\tau} \quad \text{for every real number} \,\, \alpha >0\\
			d_1,d_2 \in  \mathscr{R}_{\tau} \Rightarrow d_1 \vee d_2 \in  \mathscr{R}_{\tau} 
		\end{cases*}
	\end{equation}
	where $\vee$ means  $\max\{d_1,d_2\}$.
	Instead, when $\Pi$ is the set of all pseudo-metrics for $X$, a topological ``real" structure $\mathscr{R}_{\tau}$ is the greatest family of pseudo-metrics for $X$ defining an uniform topology $\tau$ on $X$. In this case we will say that $\mathscr{R}_{\tau}$ is a \emph{uniformizable topological  ``real" structure}. 
\end{example}

\begin{example}\label{ex:deltaloc}
	We say that  $\mathscr{P}' \prec_{\Pi} \mathscr{P}$ \emph{in the $\Delta$-local way} if $\wh{\mathscr{F}}(\mathscr{P}') \subseteq \wh{\mathscr{F}}(\mathscr{P})$. As $\mathscr{P}$ and $\mathscr{P}'$ are equivalent if they define the same $\Delta$-local filter on $X \times X$, a real structure on $X$ is, in this case, the greatest family of maps  $d : X \times X \mapsto \mathbb{R}^+$ which define the same $\Delta$-local filter on $X \times X$.
\end{example}

\begin{example}\label{ex:strogdeltaloc}
	We say that $\mathscr{P}' \prec_{\Pi} \mathscr{P}$ \emph{in the strong $\Delta$-local way} if 
	$$\wh{\mathscr{F}_x}(\mathscr{P}') \subseteq \wh{\mathscr{F}_x}(\mathscr{P}) \quad \forall x \in X.$$
	With this preorder, $\mathscr{P}$ and $\mathscr{P}'$ are equivalent if, for every $x \in X$,  $\wh{\mathscr{F}_x}(\mathscr{P}') = \wh{\mathscr{F}_x}(\mathscr{P})$. Observe that this implies $\wh{\mathscr{F}}(\mathscr{P}') = \wh{\mathscr{F}}(\mathscr{P})$, while the converse is not true.
	
	Thus, the real structures on $X$ are now the greatest families of maps  $d : X \times X \mapsto \mathbb{R}^+$ which define, for every $x \in X$, the same $\mathscr{F}_x$.
\end{example}

\begin{example}\label{ex:lipsch}
	We say that $\mathscr{P}' \prec_{\Pi} \mathscr{P}$ \emph{in the Lipschitz way} if for every $d' \in \mathscr{P}'$ there are a finite subfamily $(d_i)_{i \in I}$ of $\mathscr{P}$ and a real number $\alpha >0$ such that
	$$d' \leq \alpha \max_{i \in I} d_i.$$
	The real structure $\mathscr{R}_l$ related to this preorder will be called a \emph{weak Lipschitz structure} when $\Pi$ is the set of all maps $d : X \times X \mapsto \mathbb{R}^+$ while $\mathscr{R}_l$ will be said a \emph{Lipschitz structure} when $\Pi$ is the set of all pseudo-metrics  $d : X \times X \mapsto \mathbb{R}^+$.
	
	It is not difficult to realize that a family $\mathscr{R}_l$ of maps $d : X \times X \mapsto \mathbb{R}^+$ is a weak Lipschitz structure for $X$ if and only if
	
	\begin{equation}\label{lipstruct}
		\begin{cases*}
			d_1 \leq d, d \in \mathscr{R}_l \Rightarrow d_1 \in  \mathscr{R}_{l}\\
			d \in  \mathscr{R}_{l} \Rightarrow \alpha d \in  \mathscr{R}_{l} \quad \text{for every real number} \,\, \alpha >0\\
			d_1,d_2 \in  \mathscr{R}_{l} \Rightarrow d_1 \vee d_2 \in  \mathscr{R}_{l}. 
		\end{cases*}
	\end{equation}
	Moreover, a family $\mathscr{R}_l$ of pseudo-metrics $d : X \times X \mapsto \mathbb{R}^+$ is a Lipschitz structure for $X$ if and only if \eqref{lipstruct} is satisfied. We remark that \eqref{lipstruct} define Lipschitz structures in a similar way as in \cite{fras70}.
\end{example}

Taking into account Examples~\ref{ex:top} and \ref{ex:lipsch} we can state the following Remark.
\begin{remark}
	Each topological ``real" structure $\mathscr{R}_{\tau}$ contains a weak Lipschitz structure $\mathscr{R}_{l}$. Each uniformizable topological ``real" structure contains a Lipschitz structure. 
\end{remark}
\section{Other examples of preorders $\prec_{\Pi}$.}\label{ex2:preord}
Let $\mathscr{P}$ and $\mathscr{P}'$ be two families of pseudo-metrics on a set $X$.
\begin{example}
	We will say that \emph{$\mathscr{P}$ absorbs $\mathscr{P}'$ in an uniform way} if $\mathscr{U}(\mathscr{P}') \subseteq \mathscr{U}(\mathscr{P})$, where
	$\mathscr{U}(\mathscr{P}')$ is the uniformity generated by $\mathscr{P}'$ and $\mathscr{U}(\mathscr{P})$ is the uniformity generated by $\mathscr{P}$, namely if for every $d' \in \mathscr{P}'$ and every $\varepsilon >0$ there are a finite subfamily $(d_i)_{i \in I}$ of $\mathscr{P}$ and a number $\delta_{\varepsilon} >0$ such that 
	
	\begin{equation}\label{abs:unif}
		\bigvee_{i \in I}d_i(x_1,x_2) \leq \delta_{\varepsilon} \Rightarrow d'(x_1,x_2) \leq \varepsilon \quad \forall x_1,x_2 \in X.
	\end{equation} 
	
	The real structure defined by this preorder is the \emph{gage} of an uniformity for $X$ (see \cite[pp. 188-189]{kell55}).
\end{example}

\begin{example}
	We will say that \emph{$\mathscr{P}$ absorbs $\mathscr{P}'$ in a quasi-Lipschitz way}, if for every $d' \in \mathscr{P}'$ and every $\varepsilon >0$ there are a finite subfamily $(d_i)_{i \in I}$ of $\mathscr{P}$ and a number $\alpha_{\varepsilon} >0$ such that
	
	\begin{equation}\label{abs:quasiLip}
		d'(x_1,x_2) \leq \alpha_{\varepsilon} \bigvee_{i \in I}d_i(x_1,x_2) \quad \forall x_1,x_2 \in X.
	\end{equation} 
\end{example}

\begin{remark}
	If $\mathscr{P}$ absorbs $\mathscr{P}'$ in a quasi-Lipschitz way then $\mathscr{P}$ absorbs $\mathscr{P}'$ in an uniform way.
\end{remark}
Indeed, \eqref{abs:quasiLip} implies \eqref{abs:unif} with $\delta_{\varepsilon} =(2 \alpha_{\varepsilon})^{-1}$.

\begin{example}
	We will say that \emph{$\mathscr{P}$ absorbs $\mathscr{P}'$ in an uniform quasi-Lipschitz way} when for every $d' \in \mathscr{P}'$ and every $\varepsilon >0$ there are a finite subfamily $(d_i)_{i \in I}$ of $\mathscr{P}$ and two numbers $\alpha_{\varepsilon} >0$ and $\beta_{\varepsilon} >0$ such that
	\begin{equation}
		\bigvee_{i \in I}d_i(x_1,x_2) \leq \beta_{\varepsilon} \Rightarrow d'(x_1,x_2) \leq \alpha_{\varepsilon} \bigvee_{i \in I}d_i(x_1,x_2) + \varepsilon.
	\end{equation}
\end{example}
\begin{example}
	We will say that \emph{$\mathscr{P}$ absorbs $\mathscr{P}'$ in an local quasi-Lipschitz way} when every $x \in X$ has a neighborhood $U_x$, which depends only on $\mathscr{P}$, for the topology defined by $\mathscr{P}$, such that for every $d' \in \mathscr{P}'$ and every $\varepsilon >0$ there are 
	a finite subfamily $(d_i)_{i \in I}$ of $\mathscr{P}$ and a number $\alpha_{\varepsilon} >0$ such that
	\begin{equation}
		d'(x_1,x_2) \leq \alpha_{\varepsilon} \bigvee_{i \in I}d_i(x_1,x_2)+ \varepsilon \quad \forall x_1,x_2 \in U_x.
	\end{equation}
\end{example}
\begin{remark}
	When (as we are considering in this section) $\Pi$ is the set of the pseudo-metrics on $X$, the real structures defined by the preorder considered in 
	Example~\ref{ex:lipsch} will be called \emph{Lipschitz structures} (instead of weak Lispchitz structures).
\end{remark}

It is particularly interesting the case when $\Pi$ is the set of families of weak pseudo-metrics on $X$ (with the meaning given in the Introduction), because all topologies on $X$ and all pseudo-uniformities for $X$ are defined by some family of weak pseudo-metrics on $X$, as it is shown in \cite{val23}. It follows that, in this case,

\begin{remark}
	If $\tau$ is a topology on $X$, the real structure associated to $\tau$ is the greatest family of weak pseudo-metrics which defines $\tau$.
\end{remark}
\begin{remark}
	If $\mathscr{U}$ is a pseudo-uniformity for $X$, the real structure associated to $\mathscr{U}$ is the greatest family of weak pseudo-metrics which defines $\mathscr{U}$.
\end{remark}

\section{Some examples of morphisms between real spaces.}
Let $(X,\mathscr{R}_X)$ and $(Y,\mathscr{R}_Y)$ be two real spaces of the same type, and consider a map $f : X \mapsto Y$.
\begin{example}
	If $\mathscr{R}_X$ and $\mathscr{R}_Y$ are topological real structures, then $f$ is a morphism if and only if for every $x \in X$, $d_Y \in \mathscr{R}_Y$ and a real number $\varepsilon_{\scriptscriptstyle Y} > d_Y(f(x),f(x))$ there are $d_X \in \mathscr{R}_X$ and a real number $\varepsilon_{\scriptscriptstyle X} > d(x,x)$ such that
	$$d_X(x,\xi) \vee d_X(\xi,x) < \varepsilon_{\scriptscriptstyle X} \Rightarrow d_Y(f(x),f(\xi)) \vee d_Y(f(\xi),f(x)) < \varepsilon_{\scriptscriptstyle Y} \quad \forall \xi \in X.$$  
	
	Of course, this means that $f$ is continuous for the topologies on $X$ and $Y$ defined by $\mathscr{R}_X$ and $\mathscr{R}_Y$.
\end{example}
\begin{example}
	If $\mathscr{R}_X$ and $\mathscr{R}_Y$ are $\Delta$-local real structures (see Example~\ref{ex:deltaloc}), then $f$ is a morphism if and only if for every $d_Y \in \mathscr{R}_Y$ and a real number $\varepsilon_{\scriptscriptstyle Y} > \inf_{x \in X}d_Y(f(x),f(x))$ there are $d_X \in \mathscr{R}_X$ and a real number $\varepsilon_{\scriptscriptstyle X} > \inf_{x \in X}d_X(x,x)$ such that
	
	$$ x_1,x_2 \in X, \quad d_X(x_1,x_2) < \varepsilon_{\scriptscriptstyle X} \Rightarrow d_Y(f(x_1),f(x_2)) < \varepsilon_{\scriptscriptstyle Y}.$$
\end{example}

\begin{example}
	If $\mathscr{R}_X$ and $\mathscr{R}_Y$ are strong $\Delta$-local real structures (see Example~\ref{ex:strogdeltaloc}), then $f$ is a morphism if and only if for every $x \in X$, $d_Y \in \mathscr{R}_Y$ and a real number $\varepsilon_{\scriptscriptstyle Y} > d_Y(f(x),f(x))$ there are $d_X \in \mathscr{R}_X$ and a real number $\varepsilon_{\scriptscriptstyle X} > d_X(x,x)$ such that
	
	$$ x_1,x_2 \in X, \quad d_X(x_1,x_2) < \varepsilon_{\scriptscriptstyle X} \Rightarrow d_Y(f(x_1),f(x_2)) < \varepsilon_{\scriptscriptstyle Y}.$$
\end{example}

%\subsection*{Acknowledgements}
%Place all thanks and grant acknowledgements here.

%%%%%%%%%%% To ease editing, use normal size for the references:

\end{document}